\documentclass[twoside]{article}

\usepackage[utf8]{inputenc}
\usepackage{fancyhdr}
\usepackage[square,numbers]{natbib}
\usepackage[textwidth=6in,left=1.25in,right=1.25in,asymmetric]{geometry}
\usepackage{amsmath}
\usepackage{mathrsfs}

\usepackage{mathabx}

\usepackage{amsthm}

\usepackage{amsfonts}
\usepackage{amssymb}
\usepackage{graphicx}
\usepackage{placeins}
\usepackage{algorithm2e}
\usepackage{caption}
\usepackage{array}
\usepackage{multicol}
\usepackage{color}
\usepackage{mhequ}
\usepackage{url}
\usepackage{bbm}
\usepackage[font={small,it}]{caption}
\usepackage{stackrel}
\usepackage{enumerate}
\usepackage{datetime}
\usepackage{breqn}
\usepackage{authblk}
\usepackage{datetime}
\usepackage[colorlinks=true,citecolor=blue]{hyperref}

\theoremstyle{plain}

\theoremstyle{definition}

\newtheorem{remark}{Remark}[section]

\renewcommand{\cite}{\citet}

\newcommand{\bb}[1]{\mathbb{#1}}

\newcommand{\PR}[1]{\mathrm{pr} \left( #1 \right)}

\newcommand{\mc}[1]{\mathcal{#1}}
\newcommand{\bone}[1]{\mathbbm{1}\left\{ #1 \right\}}

\newcommand{\bigO}[1]{\mc O \left( #1 \right)}
\newcommand{\bigOc}[1]{\mc O \left\{ #1 \right\}}

\DeclareMathOperator{\argmax}{argmax}

\newcommand{\be}{\begin{equs}}
\newcommand{\ee}{\end{equs}}

\DeclareMathOperator{\pr}{pr}

\makeatletter
\providecommand*{\input@path}{}
\g@addto@macro\input@path{{./}{../input/}}
\makeatother

\graphicspath{{./}}

\title{Theoretical Limits of Record Linkage and Microclustering}

\author[1]{James E. Johndrow \thanks{johndrow@stanford.edu}}
\author[2]{Kristian Lum \thanks{kl@hrdag.org}}
\author[3]{David B. Dunson \thanks{dunson@duke.edu}}
\affil[1]{Department of Statistics, Stanford University}
\affil[2]{Human Rights Data Analysis Group}
\affil[3]{Department of Statistical Science, Duke University}
\date{\today}
\allowdisplaybreaks

\begin{document}

\maketitle

\begin{abstract}
There has been substantial recent interest in record linkage, attempting to group the records pertaining to the same entities from a large database lacking unique identifiers.  This can be viewed as a type of ``microclustering,'' with few observations per cluster and a very large number of clusters.   A variety of methods have been proposed, but there is a lack of literature providing theoretical guarantees on performance.  We show that the problem is fundamentally hard from a theoretical perspective, and even in idealized cases, accurate entity resolution is effectively impossible when the number of entities is small relative to the number of records and/or the separation among records from different entities is not extremely large.  To characterize the fundamental difficulty, we focus on entity resolution based on multivariate Gaussian mixture models, but our conclusions apply broadly and are supported by simulation studies inspired by human rights applications.  These results suggest conservatism in interpretation of the results of record linkage, support collection of additional data to more accurately disambiguate the entities, and motivate a focus on coarser inference.  For example, results from a simulation study suggest that sometimes one may obtain accurate results for population size estimation even when fine scale entity resolution is inaccurate.
\end{abstract}

\section{Introduction} \label{sec:intro}
Record linkage refers to the problem of assigning records to unique entities based on observed characteristics. One example, which is the motivating problem for this work, arises in human rights research (\cite[\S 3.3]{lum2013applications}), where there is interest in recording deaths or other human rights violations attributable to a conflict, such as the ongoing conflict in Syria.  In this setting, the data are incomplete records of violations, which usually consist of a name, date of death, and place of death. In the turbulent atmosphere accompanying a conflict, often multiple organizations record information on deaths with little communication or standardization of recording practices. To make matters worse, because these data are usually gathered from oral recollections of survivors, the data are observed with noise. Often, individuals will be known by different names to different members of their social network or community, and the surviving contact providing the data may either not know the full name and exact date and place of death of the deceased, or may incorrectly recall this information to the data gathering survey. The result is multiple databases consisting of noisy observations on features of the deceased that in some cases would not uniquely identify the individual even in the absence of noise. 

A variety of methods for record linkage have been proposed (\cite{christen2012data, winkler2006overview}). Much of the literature has focused on the theoretical framework of \cite{fellegi1969theory}. In this setup, every pair of records from two databases is compared using a discrepancy function of record features and classified as either a match, a non-match, or possibly a match. The goal is to design a decision rule that minimizes the number of possible matches for fixed match and non-match error rates. In this rubric, optimal linkage rules are functions of the ratio of the conditional likelihood of the discrepancy function given that the truth is a match or non-match. The necessity of performing pairwise comparisons leads to a combinatorial explosion in the number of records, and a related literature has focused on the construction of ``blocking rules'' to limit the number of pairwise comparisons performed (\cite{al2005blocking, bilenko2006adaptive, jaro1989advances, jaro1995probabilistic, michelson2006learning}). 

An alternative and more recent approach is to perform entity resolution through clustering, where the goal is to recover the entities (clusters) from one or more noisy observations on each entity (\citet{steorts2015bayesian, steorts2014smered, steorts2015entity, betancourt2016flexible}). Model or likelihood-based methods of this sort can be equated with mixture modeling, where the number of mixture components is large and the number of observations per component very small. Historically, the focus in mixture modeling has been on regularization that penalizes large numbers of clusters, in order to obtain a more parsimonious representation of the data generating process. Recognizing that this type of regularization is inappropriate for most record linkage problems, \cite{miller2015microclustering} define the concept of \emph{microclustering}, where the cluster sizes grow at a sublinear rate in the number of observations. They propose a Bayes nonparametric approach to clustering in this setting that takes advantage of a novel random partition process that has the microclustering property. This is applied to multinomial mixtures in \cite{betancourt2016flexible}. 

While the concept of microclustering is appropriate for most record linkage problems, there is currently a lack of literature on performance guarantees and other theoretical properties of entity resolution procedures. Because microclustering methods favor sublinear growth in cluster sizes, the number of parameters of these models can often grow at the same rate as the number of observations, so basic asymptotic properties like central limit theorems, strong laws, and consistency will fail to hold. For example, in the human rights applications that motivate \cite{miller2015microclustering}, the number of unique records per entity is thought to be very small -- generally less than 10 -- while the number of unique entities is thought to be in the thousands or hundreds of thousands. As such, it is critical to consider the finite-sample performance of microclustering in cases where the number of records per cluster is a tiny fraction of the sample size, and to obtain theoretical upper bounds on how accurate cluster-based entity resolution can possibly be when the microclustering condition holds. 

Working with simple mixture models where some of the parameters are known, we characterize the exact distributions of quantities related to entity resolution. Achievable performance is shown to be a function of entity separation and the noise level. Using these results, we provide minimal conditions for accuracy in entity resolution to be bounded away from zero asymptotically as the number of records grows. We also provide an information-theoretic bound on the best possible performance in the case where some of the entities cannot be uniquely identified from noiseless observations of the available features. These results are supported by several simulation studies. Our problem is related to the extensive literature on mixture identifiability (\cite{teicher1961identifiability, teicher1963identifiability, yakowitz1968identifiability, nobile1994bayesian, holzmann2006identifiability}) and estimation of the number of components (\cite{richardson1997bayesian, tibshirani2001estimating, day1969estimating, lo2001testing}), with the important distinction that we focus on mixtures with the ``microclustering'' property, and we are interested primarily in entity resolution, not in estimation of the parameters of the mixture.

Our results initially present a very dim view of entity resolution, and indeed, it appears that the full problem is unsolvable without further information except under very strong conditions. However, in many cases inferential interest focuses on certain summary statistics of the linked records, which may be relatively insensitive to errors in entity resolution. Motivated by the human rights application mentioned above, we consider the case where the ultimate goal of entity resolution is to recover the total number of entities in the population, given that the union of all $p$ databases being linked may not contain a complete enumeration of the entities. A variety of methods exist for this problem, which is referred to as closed population estimation, and generally use as data a relatively small contingency table that characterizes the number of unique records appearing on every possible combination of the $p$ databases. This is a $2^p$ contingency table whose entries can be denoted $n(x)$ for $x \in \{0,1\}^p$; so, for example, the entry $n(011)$ in a $2^3$ table gives the count of the number of entities that appear in the second and third databases but not the first. Naturally, the count $n(000)$ gives the number of entities not observed in any of the databases, which is unknown and the object of inference. In a simulation study, we show that relatively accurate estimation of the total population size is possible even when entity resolution is inaccurate. 

\section{Preliminaries} \label{sec:prelim}
We work primarily with Gaussian mixtures. This differs from the mixture considered in \cite{betancourt2016flexible}, which is similar to \cite{dunson2009nonparametric}, a nonparametric Bayesian model for multivariate categorical data.  Our rationale for using Gaussian mixtures is the results of \cite{johndrow2017tensor} and \cite{fienberg2007maximum}, which make clear that the maximum number of unique mixture components in the model of \cite{dunson2009nonparametric} is strictly less than $d^p$, where $d$ is the number of distinct levels of the categorical variables. Thus, it is impossible to resolve more than $d^p$ entities on the basis of $p$ categorical measurements, motivating our focus on the case of continuous outcomes that does not suffer this fundamental limitation.

Consider a mixture model with likelihood
\be
L(y \mid \theta) = \sum_{k=1}^{K} \nu_k \phi(y;\mu_k,\Sigma_k), \label{eq:like}
\ee
where $\nu \in \mc S^{K-1}$ is an element of the $K-1$ dimensional probability simplex, $\mu_k, y \in \bb R^p$, $0 \le K < \infty$, $\Sigma_k$ is a $p \times p$ positive definite matrix, and $\phi(y;\mu,\Sigma) = |2 \pi \Sigma|^{-1/2} \exp\{-(y-\mu)'\Sigma^{-1} (y-\mu)/2\}$ is the Gaussian density function. In \eqref{eq:like}, $y$ are observed entity-specific features that we will use to perform record linkage. In our motivating application, typical features are name, time/date of death, and place of death. It is natural to treat time and place as continuous variables
, and it is common to embed name into an abstract continuous space by way of a metric on text, such as Jaccard similarity or Levenshtein distance. As such, \eqref{eq:like} provides a reasonable default mixture in our setting.   

In providing an upper bound on performance in entity resolution, we focus on a case that favors good performance; in particular, we consider the task of correctly determining which mixture component generated each $y_i \sim L(y|\theta)$, for $i=1,\ldots,n$, assuming that \eqref{eq:like} is known. We focus on the estimator
\be
\widehat{k}(y) = \argmax_k \log \phi(y;\mu_k,\Sigma_k) = \argmax_k \log \phi_k(y), \label{eq:khat}
\ee
that is, we will assign $y$ to the mixture component that maximizes the likelihood. In the sequel, we study a series of cases in which the set of unknown parameters in the model is gradually expanded, which provides a set of theoretically tractable finite sample bounds on the best case performance of clustering-based approaches to entity resolution. While we focus on Gaussian mixtures for simplicity, many of the results apply equally to mixtures of any kernels that are functions of a metric on $\bb R^d$, and we make an effort to point out the obvious extensions where appropriate. 

\subsection{ An information theoretic bound }
We begin by considering the case when multiple true entities have identical values of the entity-specific parameters $(\mu_k,\Sigma_k)$. In particular, consider a situation in which we observe two complete enumerations of a population, each containing a nearly mutually exclusive set of covariates about each individual in the population. In this setting, we assume that these two lists only contain one field in common. For example, suppose one list contains each individual's name and date of birth and the other contains each individual's name, location of death, and date of death. The goal is to match each individual on the first list to the correct individual on the second list to produce a complete dataset consisting of name, date of birth, date of death, and location of death for each individual in the population.  

In locations with low entropy in the name distribution, as is the case in Syria, this list is likely to be composed of many individuals sharing exactly the same first and last name. In this section, we illustrate the limitations in the performance of record linkage when multiple entities have identical values of $(\mu_k, \Sigma_k)$ and the data are observed without noise. In the context of \eqref{eq:like}, this corresponds to the limit as the maximum eigenvalue of $\Sigma_k$ approaches zero, resulting in a mixture of delta measures. For simplicity, we focus on the case where the features are names, with the obvious parallel to the case where features are vectors in $\bb R^p$ and multiple entities have identical true values of the feature vector. 

Suppose we observe a list of names $y_i$ for $i=1,\ldots,N$, where $y_i$ takes $M < N$ unique values. Let $N_m = \sum_i \bone{y_i = \mu^*_m}$ for $m=1,\ldots,M$, where each value of $\mu^*_m$ corresponds to a unique observed name. Suppose further that $x_i \in \{1,\ldots,N_{y_i}\}$ is an \emph{unobserved} unique identifier for each person. For example, the full data could look like table \ref{tab:names} and we only observe the name column.   

\begin{table}
\centering
 \begin{tabular}{lr}
  Name ($y_i$) & identifier ($x_i$) \\
  John Smith & 1 \\
  John Smith & 2 \\
  Jane Wang & 3 \\
  Jane Wang & 4\\
  Anna Rodriguez & 5 \\
  Anna Rodriguez & 6 \\
 \end{tabular}
 \caption{Example of data for name problem} \label{tab:names}
\end{table}

The inferential goal is to assign the correct identifier to each person. We consider the case where it is known that there is exactly one record corresponding to each person, and use a random allocation procedure. In the case where multiple true entities have identical values of $\mu_k$, the estimator in \eqref{eq:khat} does not give a unique solution, since
\be
\bb P(y \mid \mu_k) = \left\{ \begin{array}{cc} 1 & y=\mu_k \\ 0 & \text{else} \end{array} \right.,
\ee
so the likelihood has identical values for all $k$ such that $\mu_k=y$.  The procedure used is to collect $I_m = \{i : y_i = \mu^*_m\}$ for each $m$ and then randomly assign the $y_i, i \in I_m$ to a permutation of $\{1,\ldots,N_m\}$ such that each $y_i$ is assigned to exactly one element of $\{1,\ldots,N_m\}$. After making this assignment, the \emph{true} value of $x_i$ is revealed and the number of correct assignments enumerated. Clearly, there are $N_m!$ ways to assign each indvidual with the same name a number between 1 and $N_m$, and only one of these assignments will be \emph{exactly} right. Let $z_m$ be the number of correct assignments with name $\mu^*_m$, and let $z = \sum_m z_m$. Then the probability of assigning every $y_i$ to its true $x_i$ is
\be
\mathrm{pr}(z = N) = \prod_{m=1}^M \frac{1}{N_m!}.
\ee
On the log scale this turns out to be very intuitive since
\be
\log\{\mathrm{pr}(z = N)\} &= \sum_{m=1}^M  \log\left( \frac{1}{N_m!} \right) \\
&= -H_Y -\sum_{m=1}^M  \frac{1}{N_m} + \bigOc{\log\left( \frac{1}{N_m} \right) }
\ee
where $H_Y$ is the entropy of the name distribution. Moreover, the distribution of $z_m$ can be described by the probability mass function
\be
\mathrm{pr}(z_m = z) = \frac{1}{N_m!} {N_m \choose z} \{!(N_m-z)\},
\ee
where for an integer $n$, $!n$ is the number of \emph{derangements} of the integers $1,\ldots,n$, that is, the number of ways to re-arrange the sequence $1,\ldots,n$ such that none of the elements of the sequence are in their original locations. We have the relation
\be
!n = \left \lfloor \frac{n!}{e} + \frac{1}{2} \right \rfloor = \text{round}(n!/e),
\ee
for $n \ge 1$, where $\lfloor \cdot \rfloor$ is the floor function; also, $!0 \equiv 1$.

We now consider the expectation of $z_m$. It turns out to be relatively easy to compute upper and lower bounds; proofs are deferred to the Appendix.
\begin{remark} \label{rem:gammabound}
\be
\frac{\Gamma(N_m,1)}{\Gamma(N_m)} \le \bb E[z_m] \le \frac{\Gamma(N_m,1)}{\Gamma(N_m)} + \frac{2^{N_m-1}}{(N_m-1)!} ,
\ee
where $\Gamma(N_m,1) = \int_{1}^{\infty} t^{N_m-1} e^{-t} dt$ is the incomplete Gamma function. 
\end{remark}

\begin{figure}
 \centering
 \includegraphics[width=0.5\textwidth]{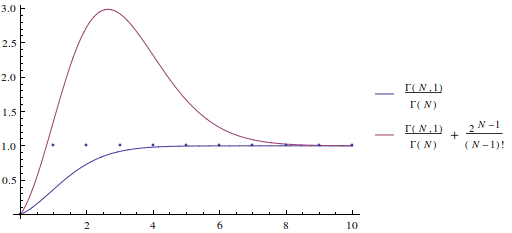}
 \caption{Upper and lower bounds on $\bb E[z_k]$ (as labeled in legend) and exact value of $\bb E[z_k]$ (points) for $N_k \in \{0,\ldots,10\}$.} \label{fig:exz}
\end{figure}

The difference between the upper and lower bounds is less than 0.001 when $N_m = 11$, so for large $N_m$ the lower bound is very accurate. Figure \ref{fig:exz} shows the upper and lower bound as well as the exact value of $\bb E[z_m]$ for $N_m \le 10$, which is identically one in all cases. From this it is clear that taking $\bb E[z_m] = 1$ for all $N_m$ is at least a very accurate approximation, and probably the exact value of the expectation. Assuming it is exact, then by independence, we have $\bb E[z] = M$, and the expected proportion correct is $M N^{-1}$. 

We now give a concentration inequality for the proportion of correct assignments, defined as the number of correct assignments $z=\sum_{m=1}^M z_m$ divided by the number of records $N$ in the list.  In particular, we have 
\be
S_M = \frac{1}{N} \sum_{m=1}^M z_m = \sum_{m=1}^{M} y_m, \quad y_m \in [0,N_m/N],
\ee
where $N = \sum_m N_m$. Then the $y_m$ are independent and $\bb E[S_M] = M/N$, so by Hoeffding's inequality
\be
\bb P[ |S_M - M/N| > t] &\le 2 \exp\left(-\frac{ -2t^2}{\sum_{m=1}^M (N_m/N)^2} \right) = 2 \exp\left(-\frac{ -2 N^2 t^2}{\sum_{m=1}^M N_m^2} \right).
\ee

We obtained data on the frequency of all surnames and given names in the United States population from the Census Bureau. Assuming independent selection of first and last names in the overall population, we estimate $\bb E[z] =0.28$ for entity resolution of the United States population on the basis of only first and last name.  Dependence between first and last names will tend to decrease this expectation.  For the United States names distribution under independence, we have $N^2 (\sum_m N_m^2)^{-1} > 6 \times 10^5$, so
\be
\bb P[ |S_M - M/N| > t] &\le 2 \exp\left(-1.2 \times 10^6 t^2 \right);
\ee
for example, the probability that $S_M > 0.29$ is less than $10^{-51}$. Hence, in the United States names example, the distribution is highly concentrated around its expectation and there is an extremely low probability of even getting one third or more of the assignments correct.  This example illustrates the fact that in many entity resolution problems, the best possible performance is strictly less than perfect accuracy due to redundancy in the true values of the entity features. This provides an upper bound on the performance achievable when features are observed with noise. 

\subsection{Analysis of noisy observations when mixture parameters are known}
Having established the limitations resulting from redundancy of the true entity features, we now analyze the effect of noise in the setting where all true entity features are distinct. We begin with a highly simplified case. Suppose we observe a data sequence $y_1,\ldots,y_N$, and that each observation originates from the mixture in \eqref{eq:like} with $K = N$, $\nu_k=N^{-1}$, $\mu_k \in [a,b]$, and $\Sigma_k = \sigma^2$ for all $k$.
Assume the parameters $\{\mu_k\}_{1 \le k \le N}, \sigma^2, \nu$ are known. On observing $y$, we make a point estimate of the mixture component that it originated from using the estimator in \eqref{eq:khat}. 
Let $k_0$ be the true value of $k$. Then  
\be
\PR{\widehat{k}(y) = k \mid k_0 = k} 
&= \pr_{\phi_k} \left\{ (y-\mu_k)^2 = \bigwedge_{\tilde{k}=1}^{N} (y-\mu_{\tilde{k}})^2  \right\}.
\ee
Now make the simplifying assumption that the $\mu_k$ are equally spaced, so that $|\mu_k - \mu_{k+1}| = \delta_N = N^{-1} |b-a| = N^{-1} \ell$ for all $k$. Then we have
\be
\PR{\widehat{k}(y) = k \mid k_0 = k} &= \pr_{\phi_k} \left\{ (y-\mu_k)^2 < \frac{\delta_N^2}{2}   \right\} = \pr_{\phi_k} \left\{ \frac{|y-\mu_k|}{\sigma} < \frac{\delta_N}{2\sigma} \right\} \\
&= \Phi\left( \frac{\delta_N}{2\sigma} \right) - \Phi\left( \frac{-\delta_N}{2\sigma} \right) = 2\Phi\left( \frac{\delta_N}{2\sigma} \right) - 1 \label{eq:trueallocprob}
\ee
for $k \ne 1, K$. For $k = 1$ or $k=K$, the expression is $\Phi(\delta_N/(2\sigma))$. When $K, N$ are large, the effect of using \eqref{eq:trueallocprob} for all $k$ is negligible, so to simplify exposition we will do so. A condition like that in \eqref{eq:trueallocprob} would hold for any mixtures where the component densities are a function of a metric on $\bb R$, with $\Phi(\cdot)$ replaced by a different distribution function. This includes many of the kernel functions commonly used in machine learning, and other common densities such as the t. 

Let $X_N$ be the number of correct classifications. For the Gaussian mixture, we have the following result. 
\begin{remark}[Infeasibility result for microclustering] \label{rem:infeasible}
 Suppose $\mu_k \in \bb R$ are equally spaced and restricted to a compact set, so that $\delta_N = N^{-1} \ell$ for $c < \infty$. Then
 \be 
 \PR{\left| \frac{X_N}{N} - \left\{ 2 \Phi \left( \frac{\ell}{2 N \sigma} \right) - 1 \right\} \right| > t } < 2 \exp(-2 t^2 N) 
 \ee
 and
 \be
 \lim_{N \to \infty} \pr(X_N = 0) = \lim_{N \to \infty} [2-2\Phi\{\ell /(2 N\sigma)\}]^N = e^{-\ell /(\sqrt{2 \pi} \sigma)}.
 \ee
\end{remark}
Therefore, in large populations, the proportion of correct assignments $N^{-1} X_N$ is highly concentrated around its expectation given by \eqref{eq:trueallocprob}, which will be very near zero when $N^{-1} \ell \ll \sigma$. Evidently, $X_N \to 0$ almost surely and the probability of zero correct assignments is bounded away from zero unless $\lim_{N\to \infty} N^{-1} \ell/\sigma>0$, which requires $\ell/\sigma = \Omega(N)$. In other words, either the width of the set containing the means must grow without bound or the observation noise must go to zero as $N$ grows. We refer to the condition $\ell/\sigma = \Omega(N)$ as \emph{infinite separation}, as it effectively requires that the entities be infinitely far apart relative to the noise level in the limit. Practically, this means that for entity resolution via microclustering, measurements on entity-specific features must get more precise as the number of entities increases. Given that this regime applies when \emph{all of the parameters of the mixture are known}, Remark \ref{rem:infeasible} suggests that the full problem of entity resolution by clustering is practically impossible in most cases. 

\subsection{The effect of dimension}
We now consider the case where dimension $p_N$ grows with $N$, and show that when the parameters of the mixture are known, infinite separation can be achieved in this regime even when the means reside on a compact set and observation noise does not decay to zero as $N \to \infty$. Consider the mixture in \eqref{eq:like} with $\mu_k \in \bb R^{p_N}$ and $\Sigma_k = \sigma^2 I_{p_N}$ for all $k$. Assume the means are restricted to the unit hypercube in $\bb R^{p_N}$ and that they are arranged so as to maximimize $\delta_N$ such that $\{\mu_k\}_{k \le N}$ is a $\delta_N$-separated set. The $\delta$-covering number of the unit hypercube $U$ in $p$ dimensions is $\mc N_{\delta}(U) = (1+1/\delta)^p$, so with $N$ points, we have $\delta_N = \Theta(N^{-1/p_N})$. Then the maximum likelihood estimator \eqref{eq:khat} satisfies
\be
\pr_{\phi_k} \left\{ \frac{(y-\mu_k)'(y-\mu_k)}{\sigma^2} < \frac{\delta_N^2}{2 \sigma^2} \right\} &\le \PR{\widehat{k}(y) = k \mid k_0=k} \le \pr_{\phi_k} \left\{ \frac{(y-\mu_k)'(y-\mu_k)}{\sigma^2} < \frac{\delta_N^2}{\sigma^2} \right\} \\
\pr\left( \chi^2_{p_N} < \frac{\delta_N^2}{2\sigma^2} \right) &\le \PR{\widehat{k}(y) = k \mid k_0=k} \le \pr\left( \chi^2_{p_N} < \frac{\delta_N^2}{\sigma^2} \right),
\ee
since the internal and external covering numbers are bounded above by $\mc N_{\delta}(U)$ and below by $\mc N_{2 \delta}(U)$ 
Appealing to the central limit theorem
\be
\pr\left( \chi^2_{p_N} < \frac{\delta_N^2}{c \sigma^2} \right) \to \Phi\left\{ \frac{\delta_N^2/(c \sigma^2) -p_N}{\sqrt{p_N/c}}  \right\},
\ee
so $\delta^2_N/\sigma^2 = \Omega(p_N)$ is a necessary and sufficient condition for $\PR{\widehat{k}(y) = k \mid k_0=k} \not{\to} 0$ as $N, p_N \to \infty$. With $N$ points, we have $\delta_N = \Theta(N^{-1/p_N})$. So $\delta_N^2/\sigma^2 = \Omega(p_N)$ implies $\sigma^2 = \Theta(p_N^{-1} N^{-2/p_N})$. Thus, we still need $\sigma^2 \to 0$, but at a much slower rate. For example, if $p_N = \log(N)$, then $\lim_{N \to \infty} N^{-2/p_N} = e^{-2}$, so it is enough to have $\sigma^2 = \{\log(N)\}^{-1}$ to achieve infinite separation. Of course, having $p_N \to \infty$ in the case where the mixture parameters are unknown means that for each mixture component, we must estimate a growing number of parameters, necessitating $K = \bigO{N p_N^{-1}}$ to even achieve consistency.  The practical ramification is that, if we ignore the need to estimate the parameters of each component, one way to combat the failure of entity resolution as the number of entities increases is to attempt to increase the number of variables collected per record on each entity.

\subsection{Case where means are unknown : Bayesian mixtures} \label{sec:bayesmix}
We now consider the case where the means are unknown but the other mixture components are known. Suppose $N$ observations are generated from the mixture
\be
y \sim \sum_{k=1}^K \pi_k \phi(\mu_k,\sigma^2),
\ee
with $\sigma^2$ and $\pi$ known. Consider a Bayesian analysis with priors $\mu_k \stackrel{iid}{\sim} \phi(0,\tau^2)$. 

Let $z_1,\ldots,z_N$ for $z_i \in \{1,\ldots,K\}$  be a configuration of the $N$ observations into $K$ classes, and let $N_k = \sum_i \bone{z_i = k}$. Let $\mc Z = \{1,\ldots,K\}^N$ be the set of all possible configurations, with $|\mc Z| = K^N$. The marginal likelihood of the configuration, integrating out the means, is
\begin{eqnarray}
 L(y,z \mid \pi, \tau^2, \sigma^2) &= \frac{N!}{\prod_k N_k!} \left( \prod_k \pi_k^{N_k} \right) \prod_{k=1}^K \frac{\sigma}{(\sqrt{2 \pi} \sigma)^{N_k} \sqrt{N_k \tau^2 + \sigma^2} } 
\nonumber \\
&\times \exp\left\{ -\frac{\sum_{i : z_i = k} y_i^2}{2 \sigma^2} + \frac{\tau^2 \left( \sum_{i : z_i = k} y_i \right)^2 }{2\sigma^2 (N_k \tau^2 + \sigma^2)} \right\},
\end{eqnarray}
so the posterior probability of the configuration is
\be
p(z \mid y, \pi, \tau^2, \sigma^2) = \frac{L(y,z \mid \pi, \tau^2, \sigma^2)}{\sum_{z^* \in \mc Z} L(y,z^* \mid \pi, \tau^2, \sigma^2)} = \frac{1}{1+\sum_{z^* \ne z} BF(z^*,z)},
\ee
where the Bayes factor is
\be
BF(z^*,z) \equiv \frac{L(y,z^* \mid \pi, \tau^2, \sigma^2)}{L(y,z \mid \pi, \tau^2, \sigma^2)}.
\ee
We now consider the case where $z^*$ consists of $N-1$ singleton clusters and a single cluster with two observations, and $z$ consists of all singletons. The Bayes factor is
\be
BF(z^*,z) = \frac{2 \pi_j}{\pi_k} \frac{\sigma \sqrt{2 \tau^2 + \sigma^2}}{(\tau^2+\sigma^2)} \exp\left[ \frac{\tau^2}{2 \sigma^2} \left\{ \frac{y_i^2+y_{i'}^2}{\tau^2+\sigma^2}- \frac{(y_i+y_{i'})^2}{2 \tau^2 + \sigma^2} \right\} \right],
\ee
where $i,i'$ are the indices of the two observations that are allocated to a single cluster in $z^*$ and two clusters in $z$, $j$ is the cluster that contains $y_i$ in configuration $z$ and is empty in configuration $z^*$, and $k$ is the index of the cluster that contains observation $y_{i'}$ in configuration $z$ and contains both $y_i,y_{i'}$ in configuration $z^*$. Integrating the Bayes factor over the data distribution, we obtain
\be
\int BF(z^*,z) \phi(y_i;\mu_i,\sigma) \phi(y_{i'}; \mu_{i'},\sigma) &= \frac{2 \pi_j}{\pi_k} \frac{2 \tau^2+\sigma^2}{\{2 (\sigma^2+\tau^2)^2 - \sigma^4\}^{1/2}} \\
&\times \exp\left[ -\frac{1}{4} \tau^2 \left\{ -\frac{(\mu_i-\mu_{i'})^2}{\sigma^4} + \frac{(\mu_i+\mu_{i'})^2}{2(\sigma^2+\tau^2)^2-\sigma^4} \right\} \right].
\ee
From this it is clear that as $||\mu_i-\mu_{i'}|| \to 0$, the expectation of the Bayes factor converges to a constant, and a necessary condition for $BF(z^*,z) \to 0$ is $||\mu_i-\mu_{i'}|| \to \infty$. Therefore, when the $\mu_i$ are confined to a compact set, Bayes factors for infinitely many incorrect configurations will converge to constants as $n \to \infty$, since $n \to \infty$ implies $||\mu_i-\mu_{i'}|| \to 0$ for infinitely many pairs $i,i'$. 

\section{Empirical analysis of entity resolution by microclustering}
We now show through simulation studies that the infeasibility results are borne out empirically. 
\subsection{Entity resolution when model parameters are known}
We first consider the case where there are $N=5,000$ entities and we observe data $y_i \sim \phi(\mu_i,\sigma^2)$ for $i=1,\ldots,N$. The common variance parameter $\sigma^2=c N^{-1}$ and $c$ is varied between $0.1$ and $2$ across the simulations, and in every case $\mu_i = i/N$, so the means are equally spaced on the unit interval. Entity resolution is performed using the estimator in \eqref{eq:khat}.  

\begin{figure}[h]
\centering
\includegraphics[width=0.5\textwidth]{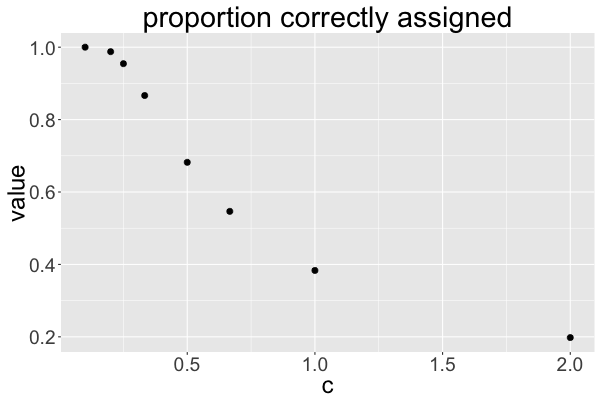}
\caption{Percent of entities correctly assigned using maximum likelihood assignment when all parameters are known. Horizontal axis shows value of $c$ used in the simulation. } \label{fig:per_right}
\end{figure}

Results are shown in Figure \ref{fig:per_right}. As expected, the proportion correctly assigned is decreasing in $c$. For $c=0.1$, entity resolution is nearly perfect, but begins to decline noticeably around $c=0.25$, which is intuitive since at that value, half the distance between the true means -- the threshold at which misassignment occurs using the maximum likelihood estimate -- is twice the standard deviation. For $c=2/3$, approximately half of the observations are correctly assigned. When $c=2$, the proportion correctly assigned is about 0.2.  

\subsection{Entity resolution when means are unknown}
\begin{figure}
 \centering
 \includegraphics[width=0.6\textwidth]{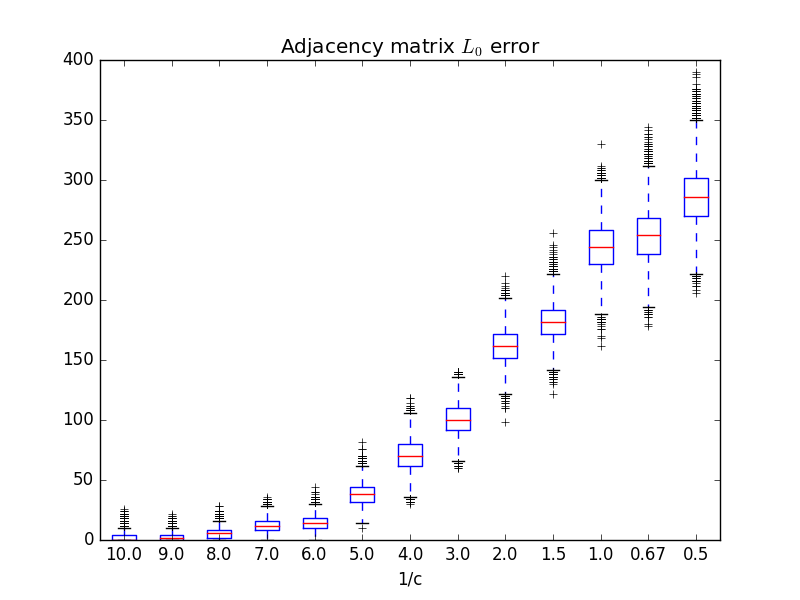}
 \caption{Boxplots of posterior samples of $||A-A^{(0)}||_0$ for Bayes mixtures with unknown means. Horizontal axis shows the reciprocal of the value of $c$ used in the corresponding simulation.} \label{fig:L0}
\end{figure}

We simulated $N=100$ observations from
\be
z_i &\sim \mathrm{Categorical}(\{1/N,\ldots,1/N\}), \quad y_i \mid z_i \sim \phi(z_i/N,\sigma^2)
\ee
for $\sigma^2=c N^{-1}$ with $c$ varied between $0.1$ and $2$ across the simulations. We then performed posterior computation by collapsed Gibbs sampling for the Bayesian mixture model with known component weights and component variances described in section \ref{sec:bayesmix}. We used identical priors $\mu_k \sim \phi(0,9)$ on the means for each component. For each MCMC sample, we compute an adjacency matrix $A$ for the 100 observations, where $A_{ij} = 1$ if observations $i,j$ are assigned to the same component, and $A_{ij}=0$ otherwise. We then compute the $L_0$ distance 
between the sampled $A$ and true adjacency matrix $A^{(0)}$
\be
||A-A^{(0)}||_0 = \sum_{i=1}^N \sum_{j=1}^N \bone{A_{ij} \ne A^{(0)}_{ij}} 
\ee
for each MCMC sample. Perfect entity resolution corresponds to $||A-A^{(0)}||_0 = 0$, while the value of $||A-A^{(0)}||_0$ can conceivably be as large as $100^2-100$, which occurs when $A$ is a matrix of ones and $A^{(0)}$ is the identity. Figure \ref{fig:L0} shows boxplots of the approximate posterior distribution of $||A-A^{(0)}||_0$ as a function of $c^{-1}$. As expected, performance in entity resolution degrades as $c$ increases, with the error rate increasing sharply near the value $c = 0.25$, as before. 

\section{Population estimation when entity resolution is poor}
In this section we give a positive empirical result. We construct a simulation in which it is possible to accurately estimate the size of a closed population from summary statistics computed on a clustering assignment even when the proportion of records correctly assigned to clusters is small.

\subsection{Simulation procedure}
Suppose each record $y_i$ is tagged with an integer $d_i \in \{1,\ldots,T\}$ for which $d_i = j$ indicates that record $i$ appeared in database $j$. Our goal is to estimate the total number of unique entities $K$ based on the entities in common in the $T$ databases. This is referred to as closed population estimation, and is the ultimate objective of record linkage in our motivating human rights setting. 

First, we describe the process for data generation. Let $K$ be the true number of entities, and for every $x \in \{0,1\}^T$, let $n(x)$ be the number of \emph{observed} records appearing in databases $j : x_j = 1$ but not in databases $j : x_j = 0$. Here $n$ is a $T$-way array of counts. We generate data from 
\be \label{eq:ctablemod}
n &\sim \text{Multinomial}(K,\pi), \quad \pi_{x} = \prod_{j=1}^T p_{j}^{x_j} (1-p_j)^{1-x_j} \\
p_j &\sim \text{Beta}(a,b), 
\ee
an independence model with a product beta prior. Then, to generate the synthetic databases, we use Algorithm \ref{alg:dbgen}.

\begin{algorithm}[H]
\SetAlgoLined
\KwData{$n(x), x \in \{0,1\}^T$}
\KwResult{$y_i, d_i$ for $i=1,\ldots, \sum_{x} n(x) \left(\sum_j x_j\right)$ }
set $\mc K = \{1,\ldots,K\}$, $N_{obs} = 0$ \;
\For{$x \in \{0,1\}^T$} { 
Sample $n(x)$ elements uniformly without replacement from $\mc K$, label result set $\mc K^*$ \;
\For{$k \in \mc K^*$} { 
$N_k  = \sum_j x_j$,  \; 
$D = \{j : x_j = 1\}$ \;
\For{$i = N_{obs}+1,\ldots,N_{obs}+N_k$}{  
generate $y_i \sim N(k/K,\sigma^2)$ \;
sample $d_i \sim \text{Discrete}(D,1)$ \;
set $D = D \setminus \{d_i\}$ \;
}
$N_{obs} = N_{obs} + N_k$ \;
}
$\mc K = \mc K \setminus \mc K^*$\;
}
\caption{Generation of synthetic databases} \label{alg:dbgen}
\end{algorithm}

This results in $T$ synthetic databases which \emph{do not contain any entries} for $n(0)$ of the $K$ entities. These are the unobserved entities, and estimating the number of unobserved entities is the goal of population estimation. In general, we choose $a,b$ in \eqref{eq:ctablemod} to make $n(0) \approx 0.25 K$.  This is consistent with real population estimation problems encountered in human rights and makes the problem relatively challenging compared to, say, the choice of $a=b=1$, which results in much smaller proportions of unobserved entities.

For the observed records $y_i, i=1,\ldots,N_{obs}$, we make a point estimate $\widehat{k}$ of the cluster assignments using \eqref{eq:khat}. For a clustering assignment $\widehat{k}_i$, we can calculate an estimate of the number of records that appear on every possible combination of the $T$ databases. For any $x \in \{0,1\}^T$, let
\be
\widehat{x}_{kj} &= \max_{i: \widehat{k}_i = k} \bone{d_i = j}, \\
\widehat{n}(x) &= \sum_k \widehat{x}_k,
\ee
giving an estimate of the list intersection counts $n(x)$ for all $x \ne (0,\ldots,0)$. We then estimate the number of unobserved entities $n(0)$ using a standard estimator for independent lists implemented in the \texttt{Rcapture} package for \texttt{R}. 

\begin{figure}
\centering
\includegraphics[width=\textwidth]{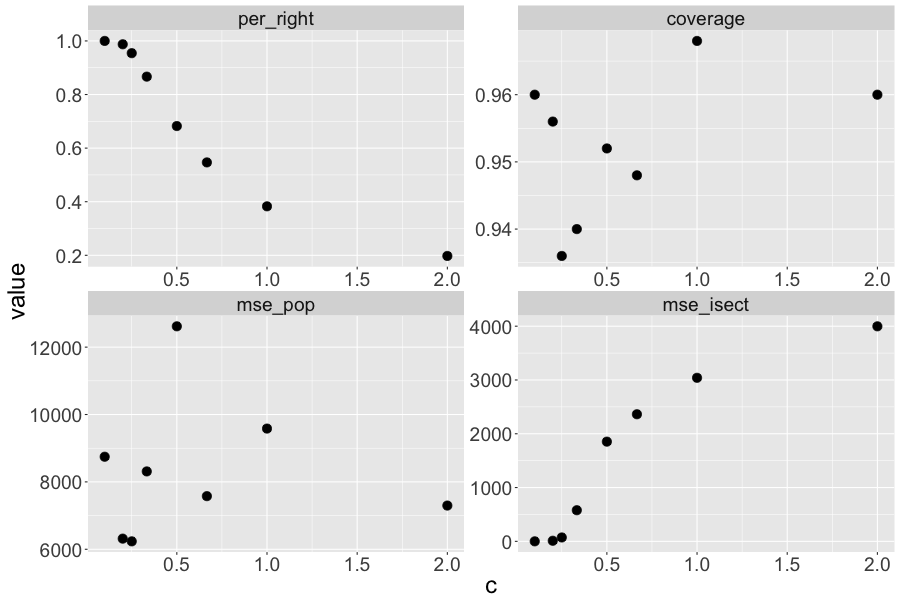}
\caption{Plots of simulation results as a function of $c$ for population estimation after entity resolution. Top left: proportion correct in entity resolution. Top right: coverage of 95 percent confidence intervals of $\widehat{n}(0)$. Bottom left: mean squared error for estimation of $n(0)$ by $\widehat{n}(0)$. Bottom right: mean squared error for estimation of $n(x)$ by $\widehat{n}(x)$ for $x \ne 0$. } \label{fig:hrsim}
\end{figure}

Results are presented in Figure \ref{fig:hrsim} for a series of simulations with $\sigma^2 = cK^{-1}$ for values of $c$ between $0.1$ and $2$. As expected, as $c$ increases, accuracy in entity resolution decreases markedly. On the other hand, coverage of $95$ percent confidence intervals for $n(0)$ and mean squared error for estimation of $n(0)$ by $\widehat{n}(0)$ are insensitive to the value of $c$. Thus, at least in this example, population estimation on the basis of linked records is not sensitive to the accuracy of entity resolution. This is particularly interesting, since estimation of $n(x)$ by $\widehat{n}(x)$ \emph{is sensitive} to the value of $c$, as shown in the bottom right panel of Figure \ref{fig:hrsim}.

\section{Discussion}

In this short article, we have shown that accurate entity resolution based on large numbers of records is essentially an unsolvable problem when the number of records per entity is very small relative to the number of entities.  This exposes a fundamental problem with entity resolution, record linkage and micro-clustering with large numbers of entities, even in idealized cases, such as when the true data-generating model is known.  These results provide a reality check on the increasing literature on entity resolution and the associated problems of record linkage and micro-clustering; it seems that any such method cannot overcome the fundamental barriers we have exposed here without some additional external information.  On the more positive side, we have also provided some simulation results illustrating that it may be possible to reliably estimate certain functionals of the linked records even if entity resolution performance is poor.  Understanding which classes of functionals we can estimate is an important area for future research.  

\section*{Acknowledgments}
This work was inspired by research conducted at Human Rights Data Analysis Group. The authors gratefully acknowledge funding support for this work from the Human Rights Data Analysis Group and the National Institutes of Health. 

\bibliographystyle{plainnat}
\bibliography{microclustering}

\begin{appendix}
\section{Appendix}
\subsection{Proof of Remark \ref{rem:gammabound}}
\be
\bb E[z_k] &\ge \sum_{j=0}^{N_k} \frac{j}{N_k!} {N_k \choose j} \frac{(N_k-j)!}{e} \\
&\ge \sum_{j=0}^{N_k} \frac{j}{j! (N_k-j)!} \frac{(N_k-j)!}{e} \\
&\ge \frac{1}{e} \sum_{j=1}^{N_k} \frac{1}{(j-1)!} = \frac{\Gamma(N_k,1)}{\Gamma(N_k)},
\ee
where $\Gamma(N_k,1) = \int_{1}^{\infty} t^{N_k-1} e^{-t} dt$ is the incomplete Gamma function. The corresponding upper bound is
\be
\bb E[z_k] &\le \sum_{j=0}^{N_k} \frac{j}{N_k!} {N_k \choose j} \left( \frac{(N_k-j)!}{e} + 1 \right) \\
&\le \sum_{j=0}^{N_k} \frac{j}{N_k!} {N_k \choose j} \left( \frac{(N_k-j)!}{e} + 1 \right) \\
&\le \frac{\Gamma(N_k,1)}{\Gamma(N_k)} + \frac{2^{N_k-1}}{(N_k-1)!}.
\ee

\subsection{Proof of Remark \ref{rem:infeasible}}
If
\be
X_N \sim \mathrm{Binomial}[N,2\Phi\{\delta_N/(2\sigma)\} - 1], \label{eq:classdist}
\ee
then
\be
\pr(X_N = 0) = [2-2\Phi\{\delta_N/(2\sigma)\}]^N. \label{eq:zeroprob}
\ee
Clearly, if the $\mu_k$ are equally spaced and restricted to be on a compact set, then $\delta_N = c/N$ for $c < \infty$. Since
\be
\lim_{N \to \infty} \pr(X_N = 0) = \lim_{N \to \infty} [2-2\Phi\{c/(2 N\sigma)\}]^N = e^{-c/(\sqrt{2 \pi} \sigma)},
\ee 
we obtain the result.

\end{appendix}
\vfill\hfill{\small Last edited: \today, \currenttime~PDT}

\end{document}